\documentclass[12pt]{article}

\usepackage{latexsym}
\usepackage{setspace}
\usepackage{amsmath}
\usepackage{amssymb}
\usepackage{latexsym}
\usepackage{amsmath}
\usepackage{amssymb}
\usepackage{natbib}
\usepackage{subfigure}

\usepackage{fancyheadings}
\setcitestyle{authoryear,open={(},close={)}}
\usepackage{graphics}
\usepackage{graphicx}
\usepackage{multirow}
\makeatletter
\def\ps@pprintTitle{%
  \let\@oddhead\@empty
  \let\@evenhead\@empty
  \let\@oddfoot\@empty
  \let\@evenfoot\@oddfoot
}
\makeatother
\makeatletter \oddsidemargin  -.5cm \evensidemargin -.5cm \textwidth
17.5cm \topmargin -.65in \textheight 24cm

\newcommand{\bd}{\begin{definition}}
\newcommand{\ed}{\end{definition}}
\newcommand{\br}{\begin{remark}}
\newcommand{\er}{\end{remark}}
\newcommand{\bea}{\begin{eqnarray}}
\newcommand{\eea}{\end{eqnarray}}
\newcommand{\beann}{\begin{eqnarray*}}
\newcommand{\eeann}{\end{eqnarray*}}
\newtheorem{theorem}{Theorem}[section]

\newtheorem{corollary}[theorem]{Corollary}
\newtheorem{remark}{Remark}[section]

\numberwithin{equation}{section}
\numberwithin{equation}{section}

\title{Estimating the scale parameters of several exponential distributions under order restriction}
\author{\small Suchandan Kayal$^{a}$, Lakshmi Kanta Patra$^{b}$\footnote
{\baselineskip=10pt
    ~patralakshmi@gmail.com}, 
 \\\small $^a$Department of Mathematics, National Institute of Technology Rourkela\\
\small $^b$Department of Mathematics, Indian Institute of Technology Bhilai}
\spacing{1.1}
\begin{document}
\date{}
\maketitle
\begin{abstract}
In the present work, we have investigated the problem of estimating parameters of several exponential distributions with ordered scale parameters under the linex loss function. We have considered estimating ordered scale parameters when the location parameters are known and unknown. For every case, we consider a class of equivariant estimators, and sufficient condition is obtained under which this class of estimators improves upon the usual estimator. Using this result, we have shown that the restricted maximum likelihood estimator is inadmissible. Finally, for every case, we conduct a simulation study to compare the risk performance of the proposed estimators. \\\\
\textbf{Keywords}: Decision theory, Best affine equivatiant estimator; Linex loss; Inadmissibility; Restricted maximum likelihood estimator.
\end{abstract}

\section{Introduction}
Researchers are often confronted with the problem of estimating parameters when they are ordered in many experimental situations, such as industrial, agriculture, clinical trials, and economic experiments. For example, it is natural to assume that the salary of a highly-skilled employee is greater than that of a low-skilled employee in a company. The average yield of a crop after using a certain type of fertilizer is expected to be more than that when no fertilizer is used. Because of various application during the last sixty years, the problem of estimating restricted parameters have drawn considerable attention from various researchers. For a good review of the problem of estimation under order restriction on the parameters, we refer to \cite{barlow1972statistical}, \cite{silvapulle2004constrained} and \cite{vaneden2006restrictd}.\\

Let $(X_{i1},\ldots, X_{in_{i}})$ be a random sample drawn from the $i$-th population $\Pi_{i},$ $i=1,\ldots,k~(k\geq2)$. The probability density function of the $i$-th population $\Pi_{i}$ is given by
\begin{eqnarray}
	f_{i}(x;\mu_{i},\sigma_{i})=\left\{\begin{array}{ll}
		\frac{1}{\sigma_{i}}~\exp{\Big(-\frac{x-\mu_{i}}{\sigma_{i}}\Big)},
		& \textrm{if $x>\mu_{i},$}\\
		0,& \textrm{otherwise,}
	\end{array} \right.
\end{eqnarray}
where $-\infty<\mu_{i}<\infty$ and $\sigma_{i}>0.$ This paper is concerned with the problem of estimating the parameters of $k~(k\geq2)$ exponential populations with ordered scale parameters. To increase system reliability, it is required to make some changes in different stages to correct the design deficiencies. Suppose there are $k$ stages in which the changes are made. Therefore, it is expected that
reliability of the system at each stage increases. Let $\lambda_{i}$ be a measure of the reliability at the $i$-th stage, $i=1,\ldots,k.$ Therefore, it is worthwhile to assume that $\lambda_{1}\leq\ldots\leq\lambda_{2}.$ Most of the work on ordered parameter estimation for exponential populations has been done with respect to the squared error loss function. For some reference, in
this direction we refer to \cite{kaur1991estimation}, \cite{vijayasree1991simultaneous,vijayasree1993mixed}, \cite{kushary1989estimating}, \cite{jin1991estimation}, \cite{pal1992order} \cite{vijayasree1995componentwise}, \cite{jana2016classification} and \cite{petropoulos2017estimation}. We will mainly focus on point estimation, and we will particularly emphasize on finding estimators which dominate
classical estimators such as maximum likelihood estimator, best affine equivariant estimator, and best scale equivariant estimator in the unrestricted problem with respect to linex loss function. The
linex loss function is given by
\begin{eqnarray}
	L(\triangle)=p'[\exp\{p\triangle\}-p\triangle-1],~\triangle=\delta-\theta,~p\neq0,~p'>0.
\end{eqnarray}
Here $p$ and $p'$ are shape and scale parameters, respectively. Throughout the paper, we assume $p'=1.$ Note that when $|p|$ tends to $0,$ then $(1.2)$ reduces to the squared error loss function. For a detailed discussion on linex loss function, we refer to \cite{zellner1986bayesian}.\\

The rest of the paper is organized as follows. In Section \ref{sec2} we have considered the estimation of ordered scale parameters when location parameters are known. We proved the inadmissibility of the best scale equivariant estimator and restricted maximum likelihood estimator. Improved estimation of ordered scale parameters with unknown location parameters is studied in Section \ref{sec3}. Finally, for every case, we have compared the relative risk performance of the proposed improved estimators. 

\section{Estimating ordered scale parameters when location parameters are known}\label{sec2}

In this section, we consider the estimation of ordered scale parameters of $k\ge 2$ exponential distributions when the location parameters are known. So without loss of generality, we assume $\mu_{i}=0,~i=1,2,\ldots,k.$
In this section, we consider estimation of $\sigma_{i}$'s subject to
the order restriction
$\sigma_{1}\leq\sigma_{2}\leq\ldots\leq\sigma_{k}.$ The loss
function is taken to be
\begin{eqnarray}
	L(\sigma_{i},\delta_{i})=\exp\left\{p\left(\frac{\delta_{i}-\sigma_{i}}
	{\sigma_{i}}\right)\right\}-p\left(\frac{\delta_{i}-\sigma_{i}}
	{\sigma_{i}}\right)-1,~p\neq0.
\end{eqnarray}
Denote $S_{i}=\sum_{j=1}^{n_{i}}X_{ij}.$ Then
$\underline{S}=(S_{1},S_{2},\ldots,S_{k})$ is a complete and
sufficient statistic for
$\underline{\sigma}=(\sigma_{1},\sigma_{2},\ldots,\sigma_{k}).$
Further, $S_{1},S_{2},\ldots,S_{k}$ are independently distributed,
where $S_{i}$ follows gamma distribution with shape parameter
$n_{i}$ and scale parameter $\sigma_{i}.$ For convenience, we denote
$S_{i}\sim Gamma(n_{i},\sigma_{i}),$ $i=1,2,\ldots,k.$ The
probability density function of $S_{i}$ is given by
\begin{eqnarray}
	f(x;n_{i},\sigma_{i})=\frac{1}{\sigma_{i}^{n_{i}}\Gamma(n_{i})}e^{-\frac{x}{\sigma_{i}}}x^{n_{i}-1},
	~x>0,
\end{eqnarray}
where $\Gamma(.)$ denotes usual gamma function. Now consider the
scale group of transformations
$G_{a_{i}}=\{g_{a_{i}}:g_{a_{i}}(x_{ij})=a_{i}x_{ij},~a_{i}>0,~j=1,2,\ldots,n_{i},
i=1,2,\ldots,k\}.$ Under this transformation, $S_{i}\rightarrow
a_{i}S_{i}$ and $\sigma_{i}\rightarrow a_{i}\sigma_{i}.$ A standard
argument shows that any scale equivariant estimator is of the form
\begin{eqnarray}
	\delta_{i_{c}}(\underline{S})=cS_{i}
\end{eqnarray}
for any constant $c.$ The following theorem provides the best scale
equivariant estimator of $\sigma_{i},$ when there is no restriction
on $\sigma_{i},$  $i=1,2,\ldots,k.$ The proof of this theorem is
omitted for the sake of brevity.
\begin{theorem}
	Under the linex loss function $(2.1),$ the (unrestricted) best scale
	equivariant estimator of $\sigma_{i}$ is
	$\delta_{i_{c_{0i}}}(\underline{S})=c_{0i}S_{i},$ where
	$c_{0i}=p^{-1}(1-e^{-p/(n_{i}+1)}),$ $i=1,2,\ldots,k.$
\end{theorem}

\subsection{Improving upon the (unrestricted) best scale equivariant estimator}

In this subsection, we take a subgroup of the scale group of
transformations $G_{a_{i}}$. Consider
$G_{a}=\{g_{a}:g_{a}(x_{ij})=ax_{ij},~a>0,~j=1,2,\ldots,n_{i},~i=1,2,\ldots,k\}.$
Then the form of the scale equivariant estimator is
\begin{eqnarray}
{\delta_{\psi}}_{i}(\underline{S})&=&S_{i}\psi_{i}\left(\frac{S_{1}}{S_{i}},\frac{S_{2}}{S_{i}},\ldots
,\frac{S_{i-1}}{S_{i}},\frac{S_{i+1}}{S_{i}},\ldots,\frac{S_{k}}{S_{i}}\right)\nonumber\\
&=&S_{i}\psi_{i}\left(\underline{U}\right),
\end{eqnarray}
where
$\underline{U}=(U_{1},U_{2},\ldots,U_{i-1},U_{i+1},\ldots,U_{k}),$
$U_{j}=S_{j}/S_{i}$ for $j\neq i,~j=1,2,\ldots,k$ and $U_{i}=S_{i}.$
The risk function of the estimators of the form $(2.4)$ is
\begin{eqnarray}
R(\sigma_{i},{\delta_{\psi}}_{i})&=&E\left[\exp\left\{\displaystyle{p\left(
\frac{S_{i}\psi_{i}(\underline{U})}{\sigma_{i}}-1\right)}\right\}-\displaystyle
p\left(
\frac{S_{i}\psi_{i}(\underline{U})}{\sigma_{i}}-1\right)-1\right]\nonumber\\
&=&E^{\underline{U}}R_{1}\left(\sigma_{i},\underline{U},\delta_{\psi_{i}}\right),
\end{eqnarray}
where
$R_{1}\left(\sigma_{i},\underline{U},{\delta_{\psi}}_{i}\right)$
represents the conditional risk function of ${\delta_{\psi}}_{i}$
given $\underline{U}=\underline{u},$ where
$\underline{u}=(u_{1},u_{2},\ldots,u_{k})$. It is given by
\begin{eqnarray}
R_{1}(\sigma_{i},\underline{U},{\delta_{\psi}}_{i})=E\left[\left(\exp\left\{\displaystyle{p\left(
\frac{S_{i}\psi_{i}(\underline{U})}{\sigma_{i}}-1\right)}\right\}-\displaystyle
p\left(
\frac{S_{i}\psi_{i}(\underline{U})}{\sigma_{i}}-1\right)-1\right)\Big|\underline{U}=\underline{u}\right].
\end{eqnarray}
Now, given $\underline{U}=\underline{u}$, the value of
$\psi_{i}(\underline{u})$ minimizing $(2.6)$ can be obtained from
the following equation:
\begin{eqnarray}
E\left[U_{i}\exp\left\{\displaystyle\frac{pU_{i}\psi_{i}(\underline{U})}{\sigma_{i}}\right\}
\Big|\underline{U} =\underline{u}\right]
=\exp\{p\}E\left(U_{i}|\underline{U}=\underline{u}\right).
\end{eqnarray}
To solve the above equation, we need to derive conditional
distribution of $U_{i}$ given $\underline{U}=\underline{u}.$ After
some calculations, it can be shown that the distribution of
$U_{i}|\underline{U}=\underline{u}$ is gamma with shape parameter
$\sum_{j=1}^{k}n_{j}$ and scale parameter
$\left(\sigma_{i}^{-1}+\sum_{j(\neq
i)=1}^{k}(u_{j}/\sigma_{j})\right)^{-1}.$ Thus, after
simplifications, from $(2.7)$ we obtain
\begin{eqnarray}
\psi_{i}(\underline{u},\underline{\sigma})=\frac{1}{p}
\left(1-\exp\left\{-\frac{p}{n+1}\right\}\right)\left(1+\sum_{j(\neq
i)=1}^{k}\frac{u_{j}\sigma_{i}}{\sigma_{j}}\right),
\end{eqnarray}
where $n=\sum_{j=1}^{k}n_{j}$. To apply \cite{brewster1974improving} 
technique, it is required to find out supremum and infimum of
$\psi_{i}(\underline{u},\underline{\sigma})$ given by $(2.8)$. For
$p\neq0,$ under the constraint
$\sigma_{1}\leq\sigma_{2}\leq\ldots\leq\sigma_{k}$, we obtain
\begin{eqnarray}
\sup_{\underline{\sigma}}\psi_{i}(\underline{u},\underline{\sigma})
&=& \left\{\begin{array}{ll}
\displaystyle\frac{1}{p}\left(1-\exp\left\{-\displaystyle\frac{p}{n+1}\right\}\right)\left(1+\sum_{j=2}^{k}
u_{j}\right),
& \textrm{if $i=1$}\\
+\infty,& \textrm{if $i=2,\ldots,k.$}
\end{array} \right.\nonumber\\
&=&{\psi_{i}}^{*}(\underline{u}),~\mbox{(say)}
\end{eqnarray}
and
\begin{eqnarray}
\inf_{\underline{\sigma}}\psi_{i}(\underline{u},\underline{\sigma})
&=& \left\{\begin{array}{ll}
\displaystyle\frac{1}{p}\left(1-\exp\left\{-\displaystyle\frac{p}{n+1}\right\}\right),
& \textrm{if $i=1$}\\
\displaystyle\frac{1}{p}\left(1-\exp\left\{-\displaystyle\frac{p}{n+1}\right\}\right)\left(1+\sum_{j=1}^{i-1}
u_{j}\right),& \textrm{if $i=2,\ldots,k$}
\end{array} \right.\nonumber\\
&=&{\psi_{i}}_{*}(\underline{u}),~\mbox{(say)}.
\end{eqnarray}
Now, we define a function $\psi_{i0}(\underline{u})$ as follows:
\begin{eqnarray}
\psi_{i0}(\underline{u}) = \left\{\begin{array}{ll}
{\psi_{i}}_{*}(\underline{u}),
& \textrm{if $\psi_{i}(\underline{u})\leq{\psi_{i}}_{*}(\underline{u})$}\\
{\psi_{i}}(\underline{u}),& \textrm{if
${\psi_{i}}_{*}(\underline{u})<\psi_{i}(\underline{u})<{\psi_{i}}^{*}(\underline{u})$}\\
{\psi_{i}}^{*}(\underline{u}),& \textrm{if
$\psi_{i}(\underline{u})\geq{\psi_{i}}^{*}(\underline{u})$}.
\end{array} \right.
\end{eqnarray}
Thus, we have the following theorem which provides sufficient
conditions for improving the estimators of the form $(2.4)$.
\begin{theorem}\label{thm2}
Let
$\delta_{\psi_{i}}(\underline{U})=S_{i}\psi_{i}\left(\underline{U}\right)$
be a scale equivariant estimator given by $(2.4)$ and $\psi_{i0}$ be
defined in $(2.11).$ If there exists some $\underline{\sigma}$ such
that
$P_{\underline{\sigma}}\left(\psi_{i0}(\underline{U})\neq\psi_{i}(\underline{U})\right)>0,$
then under the order restriction $\sigma_{1}\leq \sigma_{2}\leq
\ldots \leq \sigma_{k}$, with respect to linex loss function
$(2.1),$ the estimator $\delta_{\psi_{i0}}(\underline{U})$ improves
upon $\delta_{\psi_{i}}(\underline{U}).$
\end{theorem}
\noindent The following corollary is immediate from the Theorem
\ref{thm2}. We denote $d_{1}=p^{-1}(1-\exp\{-p/(n_{1}+n_{2}+1)\})$
and $d_{2}=p^{-1}(1-\exp\{-p/(n+1)\})$.
\begin{corollary} Suppose  it is known in prior that
$\sigma_{1}\leq \sigma_{2}\leq \ldots \leq \sigma_{k}$. Then with respect to
linex loss function given by $(2.1),$\\
$(i)$ the unrestricted best scale equivariant estimator of
$\sigma_1$ is inadmissible and dominated by the estimator given by
\begin{eqnarray}
\delta_{\psi_{10}}(\underline{U}) = \left\{\begin{array}{ll} d_{1}
S_1,
& \textrm{if $c_{01} < d_{1}$}\\
c_{01}S_1,& \textrm{if $d_{1}\le c_{01} \le
d_{1}(1+U_2)$}\\
d_{1}(1+U_2)S_1,& \textrm{if $c_{01}
>d_{1}(1+U_2)$};
\end{array} \right.
\end{eqnarray}
$(ii)$ the unrestricted best scale equivariant estimator of
$\sigma_i,$ $i=2,\dots,k$ is inadmissible and dominated by the
estimator given by
\begin{eqnarray}
\delta_{\psi_{i0}}(\underline{U}) = \left\{\begin{array}{ll}
d_{2}\left(1+\sum_{j=1}^{i-1} U_{j}\right)S_i, & \textrm{if $c_{0i}
<d_{2}\left(1+\sum_{j=1}^{i-1}
U_{j}\right)$}\\
c_{0i}S_i,& \textrm{otherwise.}
\end{array} \right.
\end{eqnarray}
\end{corollary}

\noindent Next, we consider the maximum likelihood estimator of
$\sigma_i,$ $i=1,2,\ldots,k$ under the restriction $\sigma_1 \le
\sigma_2\dots \le \sigma_k,$ which is given by
\begin{eqnarray}
\delta_{iRML}=\max_{s\le i} \min_{t \ge
i}\left\{\frac{\sum_{j=s}^tS_j}{\sum_{j=s}^tn_j}\right\}.
\end{eqnarray}
In particular, consider $k=2.$ Then the maximum likelihood
estimators of $\sigma_1$ and $\sigma_2$ under the restriction
$\sigma_1 \le \sigma_2$ reduces to
\begin{eqnarray}
\delta_{1RML}=\min\left\{\frac{Y_1}{n_1},\frac{Y_1+Y_2}{n_1+n_2}\right\}
\end{eqnarray}
and
\begin{eqnarray}
\delta_{2RML}=\max\left\{\frac{Y_2}{n_2},\frac{Y_1+Y_2}{n_1+n_2}\right\},
\end{eqnarray}
respectively. Note that the above estimators are scale equivariant
of the form $(2.4)$. Thus as an application of the Theorem $2.2$, we
get the following corollary.
\begin{corollary}
Let $\sigma_{1}\leq\sigma_{2}.$ Then for $k=2,$ under the linex loss
function $(2.1)$, the estimators
\begin{eqnarray}
\delta_{1IRML}(\underline{U}) = \left\{\begin{array}{ll} d_{1} S_1,
& \textrm{if $\phi_{01} < d_{1}$}\\
\phi_{01}S_1,& \textrm{if $d_{1}\le \phi_{01} \le  d_1 (1+U_1)$}\\
d_{2}S_1,& \textrm{if $d_{2}(1+U_{2})<\phi_{01}$};
\end{array} \right.
\end{eqnarray}
and
\begin{eqnarray}
\delta_{2IRML}(\underline{U}) = \left\{\begin{array}{ll}
d_{1}(1+U_1)S_2,
& \textrm{if $\phi_{02} < d_{1}(1+U_{1})$}\\
\phi_{02}S_2,& \textrm{otherwise,}
\end{array} \right.
\end{eqnarray}
where
$\phi_{01}=\min\left\{\frac{1}{n_1},\frac{1+U_2}{n_1+n_2}\right\}$
and
$\phi_{02}=\max\left\{\frac{1}{n_2},\frac{1+U_1}{n_1+n_2}\right\}$
have uniformly smaller risk than the restricted maximum likelihood
estimators of $\sigma_1$ and $\sigma_2$, respectively.
\end{corollary}
\subsection{A numerical study}
In this subsection, we carry out a simulation study to compare the
risk performance of the proposed estimators numerically. We have generated $50,000$ random samples from two independent exponential distributions using $R$ software for the simulation study. The
percentage relative risk improvement (PRRI) of an estimator
$\delta_{1}$ over another estimator $\delta_{2}$ is defined as
\begin{eqnarray}
	PRRI(\delta_{1},\delta_{2})=\left(\frac{Risk(\delta_{2})-Risk(\delta_1)}{Risk(\delta_{2})}\right)
	\times 100\%
\end{eqnarray}
In the following tables, we present PRRI of various estimators for
different sample sizes and different values of $p$. We have
considered sample sizes $(n_1,n_2)$ as $(3,5),$ $(10,15),$ $(8,5)$
and $(15,10).$ The values of $p$ are taken to be $0.5$ and $-0.5$.
The following points can be observed from Table $1$ and Table $2$.
\begin{enumerate}
	\item With respect to the BSEE, the PRRI of the improved BSEE
	decreases when sample sizes $(n_{1},n_{2})$ increases with the
	constraint $n_{1}<n_{2}$. Similar observation is noticed for
	$n_{1}>n_{2}$ case.
	
	\item For fixed $(n_{1},n_{2})$, the PRRI of the improved BSEE (with respect to the
	BSEE) increases if $(\sigma_{1},\sigma_{2})$ increases. Further, if
	we increase the values of $\sigma_{2}$ by keeping $\sigma_{1}$
	fixed, then the PRRI of the improved BSEE (with respect to the BSEE)
	decreases. This observation ensures that after some distance, the
	risk function of the improved BSEE coincides with that of the BSEE.
\end{enumerate}
\begin{table}[ht]
\caption{The PRRI of various estimators for $k=2$ and $p=0.5.$ The
3rd, 4th and 5th columns respectively present the PRRI of improved
BSEE, restricted MLE and improved restricted MLE over the BSEE
whereas the 6th, 7th and 8th columns present the PRRI of improved
BSEE, BSEE and improved restricted MLE over the restricted MLE.}
\vspace{.3cm}
\centering 
\scalebox{0.82}{
\begin{tabular}{c c c c c c c c} 
\hline\hline 
$(n_{1},n_{2})$ & $(\sigma_{1},\sigma_{2})$ &
$(\delta_{\psi_{10}},\delta_{\psi_{20}})$ &
$(\delta_{1RML},\delta_{2RML})$ & $(\delta_{1IRML},\delta_{2IRML})$
& $(\delta_{\psi_{10}},\delta_{\psi_{20}})$ & $(\delta_{1c_{01}},\delta_{2c_{02}})$
& $(\delta_{1IRML},\delta_{2IRML})$ \\ [0.5ex] 
\hline\hline\vspace{0.1cm} 
(3,5) & (0.2,0.5) & (3.89,6.21) & (18.26,-28.32) & (17.27,-27.01) & (-14.39,29.13)
& (-18.26,24.86) & (1.26,0.89) \\\vspace{0.1cm}
~ & (0.2,1.0) & (0.37,0.58) & (-32.19,-33.29) & (-32.43,-33.22) &
(24.63,25.41)
& (24.35,24.97) & (0.23,0.05) \\\vspace{0.1cm}
~ & (0.2,1.5) & (0.13,0.1) & (-51.64,-33.30) & (-51.84,-33.28) &
(34.14,25.05)
& (34.05,24.97) & (0.13,0.01) \\\vspace{0.1cm}
~ & (0.5,1.0) & (6.19,9.78) & (25.57,-32.73) & (24.89,-30.47) &
(-26.03,32.03)
& (-34.35,24.66) & (0.80,1.70) \\\vspace{0.1cm}
~ & (0.5,1.5) & (2.24,3.39) & (4.87,-33.21) & (3.27,-32.59) &
(-2.76,27.47)
& (-5.12,24.92) & (1.68,0.46) \\\vspace{0.1cm}
~ & (0.5,2.0) & (0.71,1.33) & (-16.48,-33.27) & (-17.39,-33.09) &
(14.76,25.97)
& (14.15,24.96) & (0.78,0.14) \\\vspace{0.1cm}
~ & (1.0,1.5) & (9.64,16.83) & (25.53,-31.53) & (29.74,-27.38) &
(-21.35,36.77)
& (-34.30,23.97) & (5.65,3.15) \\\vspace{0.1cm}
~ & (1.0,2.0) & (6.19,9.78) & (25.57,-32.73) & (24.98,-30.47) &
(-26.03,32.03)
& (-34.35,24.66) & (0.78,1.70) \\\vspace{0.1cm}
~ & (1.0,2.5) & (3.74,5.67) & (16.38,-33.08) & (15.33,-31.91) &
(-15.11,29.12) & (-19.59,24.86) & (1.25,0.88) \\ \\\\
(10,15) &
(0.2,0.5) & (0.65,0.76) & (14.7,-10.80) & (13.62,-10.61) &
(-16.48,10.44)
& (-17.24,9.75) & (1.26,0.17) \\\vspace{0.1cm}
~ & (0.2,1.0) & (0,0.002) & (-14.97,-10.80) & (-14.87,-10.81) &
(13.03,9.76)
& (13.03,9.75) & (0,0) \\\vspace{0.1cm}
~ & (0.2,1.5) & (0,0) & (-15.99,-10.81) & (-15.99,-10.80) &
(13.79,9.75)
& (13.8,9.75) & (0,0) \\\vspace{0.1cm}
~ & (0.5,1.0) & (2.53,2.90) & (27.96,-10.81) & (26.28,-9.83) &
(-35.29,12.37)
& (-38.81,9.57) & (2.32,0.88) \\\vspace{0.1cm}
~ & (0.5,1.5) & (0.13,0.19) & (0.86,-10.80) & (0.47,-10.77) &
(-0.73,9.92)
& (-0.87,9.75) & (0.38,0.034) \\\vspace{0.1cm}
~ & (0.5,2.0) & (0.02,0.015) & (-11.58,-10.81) & (-11.64,-10.80) &
(10.40,9.76)
& (10.38,9.75) & (0.051,0.003) \\\vspace{0.1cm}
~ & (1.0,1.5) & (8.66,10.72) & (13.76,-10.77) & (17.78,-6.16) &
(-5.92,19.40)
& (-15.96,9.72) & (4.56,4.16) \\\vspace{0.1cm}
~ & (1.0,2.0) & (2.53,2.91) & (27.96,-10.80) & (26.29,-9.83) &
(-35.29,12.37)
& (-38.80,9.75) & (2.32,0.88) \\\vspace{0.1cm}
~ & (1.0,2.5) & (0.65,0.76) & (14.71,-10.80) & (13.62,-10.61) &
(-16.48,10.44) & (-17.24,9.75) & (1.26,0.17) \\ \\\\
(8,5) & (0.2,0.5)
& (2.36,6.41) & (-13.92,-33.35) & (9.52,-32.10) & (14.29,29.87)
& (12.22,25.06) & (20.58,1.01) \\\vspace{0.1cm}
~ & (0.2,1.0) & (0.167,0.42) & (-0.23,-34.23) & (2.01,-34.16) &
(0.40,25.82)
& (0.23,25.5) & (2.22,0.05) \\\vspace{0.1cm}
~ & (0.2,1.5) & (0.03,0.05) & (-11.55,-34.25) & (-11.14,-34.24) &
(10.38,25.55)
& (10.35,25.51) & (0.36,0.01) \\\vspace{0.1cm}
~ & (0.5,1.0) & (4.16,12.74) & (-47.64,-32.24) & (-4.29,-29.51) &
(35.08,34.02)
& (32.27,24.38) & (29.36,2.06) \\\vspace{0.1cm}
~ & (0.5,1.5) & (1.37,3.38) & (0.30,-33.92) & (13.65,-33.21) &
(1.07,27.85)
& (-0.30,25.33) & (13.38,0.52) \\\vspace{0.1cm}
~ & (0.5,2.0) & (0.51,1.10) & (4.32,-34.18) & (9.39,-33.97) &
(-3.98,26.30)
& (-4.51,25.47) & (5.30,0.16) \\\vspace{0.1cm}
~ & (1.0,1.5) & (7.07,25.77) & (-116.0,-28.55) & (-31.36,-23.05) &
(56.98,42.26)
& (53.70,22.21) & (39.18,4.28) \\\vspace{0.1cm}
~ & (1.0,2.0) & (4.16,12.74) & (-47.64,-32.24) & (-4.29,-29.51) &
(35.08,34.02)
& (32.27,24.38) & (29.36,2.06) \\\vspace{0.1cm}
~ & (1.0,2.5) & (2.36,6.41) & (-13.92,-33.45) & (9.52,-32.10) &
(14.29,29.87) & (12.22,25.06) & (20.58,1.01) \\ \\\\
(15,10) &
(0.2,0.5) & (0.66,1.47) & (-5.75,-15.33) & (5.73,-14.91) &
(6.01,14.57)
& (5.44,13.29) & (10.86,0.36) \\\vspace{0.1cm}
~ & (0.2,1.0) & (0,0.001) & (-2.63,-15.33) & (-2.56,-15.33) &
(2.56,13.29) & (2.56,13.29) & (0.065,0) \\\vspace{0.1cm}
~ & (0.2,1.5) & (0,0) & (-9.51,-15.33) & (-9.51,-15.3) &
(9.05,13.29) & (9.05,13.29) & (0,0) \\\vspace{0.1cm}
~ & (0.5,1.0) & (2.13,5.05) & (-61.69,-15.23) & (-28.66,-13.47) &
(39.47,17.60) & (38.15,13.22) & (20.42,1.53) \\\vspace{0.1cm}
~ & (0.5,1.5) & (0.17,0.43) & (10.06,-15.33) & (14.22,-15.25) &
(-11.0,13.67) & (-11.18,13.29) & (4.62,0.072) \\\vspace{0.1cm}
~ & (0.5,2.0) & (0.004,0.03) & (5.62,-15.34) & (6.18,-15.33) &
(-5.95,13.32) & (-5.96,13.29) & (0.59,0.002) \\\vspace{0.1cm}
~ & (1.0,1.5) & (6.38,17.11) & (-195.63,-14.12) & (-99.18,-7.55) &
(68.33,27.37) & (66.17,12.37) & (32.62,5.75) \\\vspace{0.1cm}
~ & (1.0,2.0) & (2.13,5.05) & (-61.69,-15.23) & (-28.66,-13.47) &
(39.47,17.61) & (38.15,13.22) & (20.42,1.52) \\\vspace{0.1cm}
~ & (1.0,2.5) & (0.66,1.47) & (-5.75,-15.32) & (5.73,-14.91) &
(6.06,14.57) & (5.44,13.29) & (10.86,0.36) \\
[1ex] 
\hline\hline 
\end{tabular}}
\label{tb1} 
\end{table}
\clearpage

\begin{table}[ht]
\caption{The PRRI of various estimators for $k=2$ and $p=-0.5.$ The
3rd, 4th and 5th columns respectively present the PRRI of improved
BSEE, restricted MLE and improved restricted MLE over the BSEE
whereas the 6th, 7th and 8th columns present the PRRI of improved
BSEE, BSEE and improved restricted MLE over the restricted MLE.}
\vspace{.3cm}
\centering 
\scalebox{0.82}{
\begin{tabular}{c c c c c c c c} 
\hline\hline 
$(n_{1},n_{2})$ & $(\sigma_{1},\sigma_{2})$ &
$(\delta_{\psi_{10}},\delta_{\psi_{20}})$ &
$(\delta_{1RML},\delta_{2RML})$ & $(\delta_{1IRML},\delta_{2IRML})$
& $(\delta_{\psi_{10}},\delta_{\psi_{20}})$ &
$(\delta_{1c_{01}},\delta_{2c_{02}})$
& $(\delta_{1IRML},\delta_{2IRML})$ \\ [0.5ex] 
\hline\hline\vspace{0.1cm} 
(3,5) & (0.2,0.5) & (4.49,5.92) & (20.66,-10.19) & (18.62,-7.99) &
(-20.38,14.63)
& (-26.04,9.25) & (2.56,2.00) \\\vspace{0.1cm}
~ & (0.2,1.0) & (0.46,0.61) & (-4.43,-10.47) & (-4.76,-10.31) &
(4.68,10.04)
& (4.24,9.48) & (0.31,0.15) \\\vspace{0.1cm}
~ & (0.2,1.5) & (0.13,0.11) & (-13.15,-10.48) & (-13.35,-10.45) &
(11.74,9.58)
& (11.62,9.48) & (0.17,0.02) \\\vspace{0.1cm}
~ & (0.5,1.0) & (7.28,10.13) & (24.31,-9.73) & (24.13,-5.78) &
(-22.49,18.10)
& (-32.12,8.87) & (0.24,3.59) \\\vspace{0.1cm}
~ & (0.5,1.5) & (2.79,3.57) & (14,66,-10.35) & (12.77,-9.11) &
(-13.9,12.62) & (-17.18,9.38) & (2.21,1.13) \\\vspace{0.1cm}
~ & (0.5,2.0) & (1.07,1.42) & (3.39,-10.44) & (2.02,-10.03) &
(-2.40,10.74)
& (-3.51,9.46) & (1.42,0.37) \\\vspace{0.1cm}
~ & (1.0,1.5) & (11.19,17.19) & (18.68,-8.24) & (25.19,1.50) &
(-9.21,23.39)
& (-22.97,7.61) & (8.0,6.22) \\\vspace{0.1cm}
~ & (1.0,2.0) & (7.28,10.12) & (24.31,-9.73) & (24.13,-5.78) &
(-22.49,18.10)
& (-32.12,8.17) & (0.23,3.59) \\\vspace{0.1cm}
~ & (1.0,2.5) & (4.49,5.92) & (20.66,-10.19) & (18.62,-7.99) &
(-20.38,14.63) & (-26.40,9.25) & (2.56,2.00) \\ \\\\
(10,15) & (0.2,0.5) & (0.75,0.77) & (16.79,-3.68) & (15.28,-3.35) &
(-19.28,4.29)
& (-20.18,3.55) & (1.81,0.32) \\\vspace{0.1cm}
~ & (0.2,1.0) & (0,0.002) & (-4.67,-3.68) & (-4.67,-3.68) &
(4.63,3.55)
& (4.63,3.55) & (0.0) \\\vspace{0.1cm}
~ & (0.2,1.5) & (0,0) & (-5.32,-3.68) & (-5.32,-3.68) & (5.05,3.55)
& (5.05,3.55) & (0,0) \\\vspace{0.1cm}
~ & (0.5,1.0) & (2.8,2.92) & (25.65,-3.68) & (24.12,-2.17) &
(-30.74,6.37)
& (-34.51,3.55) & (2.06,1.46) \\\vspace{0.1cm}
~ & (0.5,1.5) & (0.18,0.19) & (6.65,-3.68) & (5.99,-3.61) &
(-6.39,3.73)
& (-7.12,3.55) & (0.70,0.072) \\\vspace{0.1cm}
~ & (0.5,2.0) & (0.02,0.02) & (-2.33,-3.69) & (-2.38,-3.68) &
(2.30,3.57)
& (2.28,3.55) & (0.05,0.005) \\\vspace{0.1cm}
~ & (1.0,1.5) & (9.67,10.71) & (5.51,-3.64) & (12.33,2.77) &
(4.41,13.85)
& (-5.83,3.51) & (7.22,6.19) \\\vspace{0.1cm}
~ & (1.0,2.0) & (2.79,2.92) & (25.65,-3.68) & (24.12,-2.17) &
(-30.75,6.37)
& (-34.51,3.55) & (2.07,1.46) \\\vspace{0.1cm}
~ & (1.0,2.5) & (0.75,0.77) & (16.79,-3.68) & (15.27,-3.35) &
(-19.28,4.29) & (-20.18,3.55) & (1.81,0.32) \\ \\\\
(8,5) & (0.2,0.5) & (2.66,6.57) & (-26.18,-9.71) & (4.18,-7.42) &
(22.85,14.84)
& (20.75,8.84) & (24.53,2.07) \\\vspace{0.1cm}
~ & (0.2,1.0) & (0.22,0.44) & (4.84,-10.74) & (7.76,-10.61) &
(-4.85,10.09)
& (-5.08,9.69) & (3.07,0.12) \\\vspace{0.1cm}
~ & (0.2,1.5) & (0.04,0.05) & (-1.30,-10.77) & (-0.77,-10.75) &
(1.32,9.77)
& (1.29,9.72) & (0.52,0.02) \\\vspace{0.1cm}
~ & (0.5,1.0) & (4.66,12.96) & (-69.06,-8.16) & (-11.62,-3.79) &
(43.60,19.53)
& (40.85,7.54) & (33.97,0.04) \\\vspace{0.1cm}
~ & (0.5,1.5) & (1.55,3.53) & (-6.13,-10.32) & (11.44,-9.08) &
(7.25,12.55)
& (5.78,9.35) & (16.56,1.12) \\\vspace{0.1cm}
~ & (0.5,2.0) & (0.59,1.16) & (5.17,-10.67) & (11.82,-10.27) &
(-4.83,10.69)
& (-5.44,9.64) & (7.01,0.67) \\\vspace{0.1cm}
~ & (1.0,1.5) & (7.87,25.78) & (-105.43,-3.67) & (-41.70,4.48) &
(63.79,28.42)
& (60.69,4.54) & (44.30,7.86) \\\vspace{0.1cm}
~ & (1.0,2.0) & (4.65,12.96) & (-69.07,-8.17) & (-11.62,-3.79) &
(43.61,19.53)
& (40.85,7.54) & (33.97,4.04) \\\vspace{0.1cm}
~ & (1.0,2.5) & (2.65,6.57) & (-26.18,-9.71) & (4.78,-7.43) &
(22.86,14.84) & (20.74,8.84) & (24.54,2.07) \\ \\\\
(15,10) & (0.2,0.5) & (0.73,1.47) & (-15.12,-4.97) & (-0.25,-4.29) &
(13.77,6.14) &
 (13.13,4.74) & (12.91,0.05) \\\vspace{0.1cm}
~ & (0.2,1.0) & (0,0) & (2.41,-4.99) & (2.49,-4.99) &
(-2.47,4.75) & (-2.47,4.75) & (0.08,0) \\\vspace{0.1cm}
~ & (0.2,1.5) & (0,0) & (-3.02,-4.99) & (-3.02,-4.99) &
(2.94,4.75) & (2.94,4.75) & (0,0) \\\vspace{0.1cm}
~ & (0.5,1.0) & (2.31,5.00) & (-82.77,-4.85) & (-39.81,-2.28) &
(46.54,9.41) & (45.28,4.63) & (23.50,2.45) \\\vspace{0.1cm}
~ & (0.5,1.5) & (0.02,0.40) & (6.50,-4.99) & (11.87,-4.84) &
(-6.74,5.13) & (-6.94,4.75) & (5.74,0.14) \\\vspace{0.1cm}
~ & (0.5,2.0) & (0.007,0.03) & (7.99,-4.99) & (8.72,-4.98) &
(-8.68,4.77) & (-8.69,4.75) & (0.79,0.01) \\\vspace{0.1cm}
~ & (1.0,1.5) & (6.83,16.87) & (-224.95,-3.46) & (-118.74,5.39) &
(72.99,19.66) & (71.01,3.35) & (36.58,8.56) \\\vspace{0.1cm}
~ & (1.0,2.0) & (2.30,5.00) & (-82.77,-4.85) & (-39.81,-2.82) &
(46.54,9.41) & (45.28,4.63) & (23.50,2.45) \\\vspace{0.1cm}
~ & (1.0,2.5) & (0.73,1.47) & (-15.12,-4.97) & (-0.25,-4.29) &
(13.77,6.14) & (13.13,4.74) & (12.91,0.64) \\
[1ex] 
\hline\hline 
\end{tabular}}
\label{tb1} 
\end{table}
\clearpage

\section{Estimation of ordered scale parameters when location parameters are unknown}\label{sec3}
 Throughout the present
section, we assume that $\sigma_{i}$'s and $\mu_{i}$'s,
$i=1,\ldots,k$ are unknown. Further, we assume an additional
condition on the scale parameters as
$\sigma_{1}\leq\ldots\leq\sigma_{k}$. Let us denote
$$X_{i(1)}=\min_{1\leq j\leq
n_{i}}X_{ij},~~T_{i}=\sum_{j=1}^{n_{i}}\left(X_{ij}-X_{i(1)}\right),$$
where $i=1,\ldots,k.$ It is known that
$\left(\underline{X_{(1)}},\underline{T}\right)$ is a complete and
sufficient statistic, where
$\underline{X_{(1)}}=(X_{1(1)},\ldots,X_{k(1)})$ and
$\underline{T}=(T_{1},\ldots,T_{k})$. Moreover,
$(X_{1(1)},\ldots,X_{k(1)},T_{1},\ldots,T_{k})$ are independently
distributed with $X_{i(1)}\sim
E\left(\mu_{i},\frac{\sigma_{i}}{n_{i}}\right)$ and $T_{i}\sim
Gamma(n_{i}-1,\sigma_{i})$, where $i=1,\ldots,k.$

We consider the problem of estimation of the
scale parameters $\sigma_{i}$ associated with the order restriction
$\sigma_{1}\leq\sigma_{2}\leq\ldots\leq\sigma_{k}$ when the location
parameters $\mu_{i}$'s are unknown. In this direction, we take linex
loss function given by $(2.1)$. Consider the affine group of
transformations
$G_{a_{i},b_{i}}=\{g_{a_{i},b_{i}}(x)=a_{i}x_{ij}+b_{i},~j=1,\ldots,n_{i},~i=1,\ldots,k\}$.
Now under this group, the form of the affine equivariant estimators
can be obtained as
\begin{eqnarray}
\delta_{id}(\underline{X_{(1)}},\underline{T})=dT_{i},
\end{eqnarray}
where $d$ is a constant. The following theorem provides the best
affine equivariant estimator of $\sigma_{i}$ when there is no
restriction on $\sigma_{i}$'s.
\begin{theorem}
Under linex loss function given by $(2.1),$ the (unrestricted) best
affine equivariant estimator of $\sigma_{i}$ is
$\delta_{id_{0i}}(\underline{X_{(1)}},\underline{T})=d_{0i}T_{i}$,
where $d_{0i}=\frac{1}{p}(1-e^{-p/n_{i}})~i=1,\ldots,k.$
\end{theorem}
\noindent{\bf Proof.} The proof is simple and hence omitted for the
sake of brevity. \hfill\(\Box\)

\subsection{Improving upon the (unrestricted) best affine
equivariant estimator}

Consider a subgroup of $G_{a_{i},b_{i}}$ as $G_{a,b_{i}}$. Under the
transformation
$g_{a,b_{i}}=ax_{ij}+b_{i},~i=1,\ldots,k,~j=1,\ldots,n_{i},$ the
form of the estimators is
\begin{eqnarray}
\delta_{\phi_{i}}\left(\underline{X_{(1)}},\underline{T}\right)
&=&T_{i}\phi_{i}\left(\frac{T_{1}}{T_{i}},
\ldots,\frac{T_{i-1}}{T_{i}},\frac{T_{i+1}}{T_{i}},\ldots,\frac{T_{k}}{T_{i}}\right)\nonumber\\
&=&T_{i}\phi_{i}\left(\underline{V}\right),
\end{eqnarray}
where $\underline{V}=(V_{1},\ldots,V_{i-1},V_{i+1},\ldots,V_{k}),$
$V_{j}=T_{j}/T_{i}$ for $j\neq i$ and $V_{i}=T_{i}.$ The risk of the
estimators of the form $(3.2)$ is
\begin{eqnarray}
R\left(\sigma_{i},{\delta_{\phi}}_{i}\right)&=&E\left[\exp\Big\{p\Big(
\frac{T_{i}\phi_{i}\left(\underline{V}\right)}{\sigma_{i}}-1\Big)\Big\}-\displaystyle
p\Big(
\frac{T_{i}\phi_{i}(\underline{V})}{\sigma_{i}}-1\Big)-1\right]\nonumber\\
&=&E^{\underline{V}}R_{1}\left(\sigma_{i},\underline{V},{\delta_{\phi_{i}}}\right),
\end{eqnarray}
where
$R_{1}\left(\sigma_{i},\underline{V},{\delta_{\phi}}_{i}\right)$
represents the conditional risk of ${\delta_{\phi}}_{i}$ given
$\underline{V}=\underline{v}$. It is given by
\begin{eqnarray}
R_{1}\left(\sigma_{i},\underline{V},{\delta_{\phi}}_{i}\right)=E\left[\left(\exp\left\{p\left(
\frac{T_{i}\phi_{i}\left(\underline{V}\right)}{\sigma_{i}}-1\right)\right\}-\displaystyle
p\left(
\frac{T_{i}\phi_{i}(\underline{V})}{\sigma_{i}}-1\right)-1\right)\Big|\underline{V}=\underline{v}\right].
\end{eqnarray}
Note that $R_{1}(.,.,.)$ is a convex function in $\phi_{i}(.)$.
Therefore, the choice of $\phi_{i}(.)$ which minimizes
$R_{1}(.,.,.)$ given by $(3.4)$ can be obtained from the following
relation:
\begin{eqnarray}
E\left[T_{i}\exp\left\{\left(
\frac{pT_{i}\phi_{i}\left(\underline{V}\right)}{\sigma_{i}}\right)\right\}\Big|\underline{V}
=\underline{v}\right]
=\exp\{p\}E\left(T_{i}|\underline{V}=\underline{v}\right).
\end{eqnarray}
To get $\phi_{i}\left(\underline{v}\right)$ from $(3.5)$, it is
required to derive the distribution of
$T_{i}|\underline{V}=\underline{v},~i=1,\ldots,k.$ It can be shown
that $T_{i}|\underline{V}=\underline{v}\sim
Gamma\left(n-k,(\frac{1}{\sigma_{i}}+\sum_{j(\neq
i)=1}^{k}\frac{v_{j}}{\sigma_{j}})^{-1}\right)$. Thus after some
simplification, from $(3.5)$ we obtain
\begin{eqnarray}
\phi_{i}\left(\underline{v},\underline{\sigma}\right)=\frac{1}{p}
\left(1-\exp\left\{-\displaystyle\frac{p}{n-k+1}\right\}\right)\left
(1+\sum_{j(\neq i)=1}^{k}\frac{v_{j}\sigma_{i}}{\sigma_{j}}\right).
\end{eqnarray}
Next, we obtain the infimum and supremum of
$\phi_{i}\left(\underline{v},\underline{\sigma}\right)$ given by
$(3.6)$ with respect to $\underline{\sigma}$. Under the order
restriction $\sigma_{1}\leq\ldots\leq\sigma_{k}$ and $p\neq0,$ the
infimum and supremum of
$\phi_{i}\left(\underline{v},\underline{\sigma}\right)$ can be
obtained as
\begin{eqnarray}
\inf_{\underline{\sigma}}\phi_{i}\left(\underline{v},\underline{\sigma}\right)
&=& \left\{\begin{array}{ll}
\displaystyle\frac{1}{p}\left(1-\exp\left\{-\displaystyle\frac{p}
{n-k+1}\right\}\right),
& \textrm{if $i=1$}\\
\displaystyle\frac{1}{p}\left(1-\exp\left\{-\displaystyle\frac{p}{n-k+1}\right\}\right)
\left(1+\sum_{j=1}^{i-1} v_{j}\right),& \textrm{if $i=2,\ldots,k$}
\end{array} \right.\nonumber\\
&=&{\phi_{i}}_{*}\left(\underline{v}\right),~\mbox{(say)}
\end{eqnarray}
and
\begin{eqnarray}
\sup_{\underline{\sigma}}\phi_{i}\left(\underline{v},\underline{\sigma}\right)
&=& \left\{\begin{array}{ll}
\displaystyle\frac{1}{p}\left(1-\exp\left\{-\displaystyle\frac{p}{n-k+1}\right\}\right)
\left(1+\sum_{j=2}^{k} v_{j}\right),
& \textrm{if $i=1$}\\
\infty,& \textrm{if $i=2,\ldots,k$}
\end{array} \right.\nonumber\\
&=&{\phi_{i}}^{*}\left(\underline{v}\right),~\mbox{(say)}.
\end{eqnarray}
Now, we define the following function
\begin{eqnarray}
\phi_{i0}\left(\underline{v}\right) = \left\{\begin{array}{ll}
{\phi_{i}}_{*}\left(\underline{v}\right),
& \textrm{if $\phi_{i}\left(\underline{v}\right)\leq{\phi_{i}}_{*}\left(\underline{v}\right)$}\\
{\phi_{i}}\left(\underline{v}\right),& \textrm{if
${\phi_{i}}_{*}\left(\underline{v}\right)<\phi_{i}\left(\underline{v}\right)
<{\phi_{i}}^{*}\left(\underline{v}\right)$}\\
{\phi_{i}}^{*}\left(\underline{v}\right),& \textrm{if
$\phi_{i}\left(\underline{v}\right)\geq{\phi_{i}}^{*}\left(\underline{v}\right)$}.
\end{array} \right.
\end{eqnarray}
Thus we have the following result which provides the sufficient
conditions for the improvement of the estimators of the form given
by $(3.2)$.
\begin{theorem}
Consider the estimator of the form
$\delta_{\phi_{i}}\left(\underline{X_{(1)}},\underline{T}\right)=T_{i}\phi_{i}(\underline{V})$.
Let $\phi_{i0}(\underline{v})$ be a function as defined in $(3.9)$.
If there exists some $\underline{\sigma}$ such that
$P_{\underline{\sigma}}\left(\phi_{i0}\left(\underline{V}\right)
\neq\phi_{i}\left(\underline{V}\right)\right)>0,$ then under the
order restriction $\sigma_{1}\leq \sigma_{2}\leq \ldots \leq
\sigma_{k}$,  the estimator
$\delta_{\phi_{i0}}\left(\underline{X_{(1)}},\underline{T}\right)$
dominates
$\delta_{\phi_{i}}\left(\underline{X_{(1)}},\underline{T}\right)$
with respect to the linex loss function $(2.1).$
\end{theorem}

In the following we present a corollary which is immediately follows from the
Theorem $3.2.$ For convenience of the presentation, we denote
$e_{1}=p^{-1}\left(1-\exp\{-p/(n_{1}+n_{2}-1)\}\right)$ and
$e_{2}=p^{-1}\left(1-\exp\{-p/(n-k+1)\}\right)$.
\begin{corollary}
Let $\sigma_{1}\leq\sigma_{2}\leq\ldots\leq\sigma_{k}$. Then under
the loss function given by $(2.1)$ \\
$(i)$ the unrestricted best affine equivariant estimator of
$\sigma_1$ is inadmissible and dominated by the estimator given by
\begin{eqnarray}
\delta_{\phi_{10}}\left(\underline{X_{(1)}},\underline{T}\right) =
\left\{\begin{array}{ll} e_{1} T_1,
& \textrm{if $d_{01} < e_{1}$}\\
d_{01}T_1,& \textrm{if $e_{1}\le d_{01} \le
e_{1}(1+V_2)$}\\
e_{1}(1+V_2)T_1,& \textrm{if $d_{01}
>e_{1}(1+V_2)$};
\end{array} \right.
\end{eqnarray}
$(ii)$ the unrestricted best affine equivariant estimator of
$\sigma_i,$ $i=2,\dots,k$ is inadmissible and dominated by the
estimator given by
\begin{eqnarray}
\delta_{\phi_{i0}}\left(\underline{X_{(1)}},\underline{T}\right) =
\left\{\begin{array}{ll} e_{2}\left(1+\sum_{j=1}^{i-1}
V_{j}\right)T_i, & \textrm{if $d_{0i} <e_{2}\left(1+\sum_{j=1}^{i-1}
V_{j}\right)$}\\
d_{0i}T_i,& \textrm{otherwise.}
\end{array} \right.
\end{eqnarray}
\end{corollary}
Now we consider the maximum
likelihood estimation of $\sigma_{i},~i=1,\ldots,k$ when there is an
order restriction on $\sigma_{i}$'s. Under the restriction
$\sigma\leq \sigma_{2}\leq\ldots\leq \sigma_{k}$, the restricted
maximum likelihood estimator is given by
\begin{eqnarray}
\delta^{2}_{iRML}=\max_{s\leq i}\min_{t\geq
i}\left\{\frac{\sum_{j=s}^{t}T_{j}}{\sum_{j=s}^{t}n_{j}}\right\},~i=1,\ldots,k.
\end{eqnarray}
Let $k=2.$ Then the the restricted maximum likelihood estimators of
$\sigma_{1}$ and $\sigma_{2}$, while $\sigma_{1}\leq\sigma_{2}$ are
respectively given by
\begin{eqnarray}
\delta^{2}_{1RML}=\min\left\{\frac{T_{1}}{n_{1}},\frac{T_{1}+T_{2}}{n_{1}+n_{2}}\right\}
\end{eqnarray}
and
\begin{eqnarray}
\delta^{2}_{2RML}=\max\left\{\frac{T_{2}}{n_{2}},\frac{T_{1}+T_{2}}{n_{1}+n_{2}}\right\}.
\end{eqnarray}
Note that the estimators given by $(3.13)$ and $(3.14)$ belongs to
the class of affine equivariant estimators $(3.1)$. Thus as an
application of the Theorem $3.2,$ one can easily obtain
corresponding improved estimators in terms of the risk function as
follows.
\begin{corollary}
Assume that $\sigma_{1}\leq\sigma_{2}.$ Then for $k=2,$ under the
linex loss function $(2.1),$ the estimators
\begin{eqnarray}
\delta^{2}_{1IRML}\left(\underline{X_{(1)}},\underline{T}\right) =
\left\{\begin{array}{ll} e_{1} T_1,
& \textrm{if $\xi_{01}(\underline{V}) < e_{1}$}\\
\xi_{01}(\underline{V})T_1,& \textrm{if $e_{1}\le \xi_{01}(\underline{V}) \le  e_1 (1+V_2)$}\\
e_{1}(1+V_{2})T_1,& \textrm{if
$e_{1}(1+V_{2})<\xi_{01}(\underline{V})$};
\end{array} \right.
\end{eqnarray}
and
\begin{eqnarray}
\delta^{2}_{2IRML}\left(\underline{X_{(1)}},\underline{T}\right) =
\left\{\begin{array}{ll} e_{1}(1+V_1)T_2,
& \textrm{if $\xi_{02}(\underline{V}) < e_{1}(1+V_{1})$}\\
\xi_{02}(\underline{V})T_2,& \textrm{otherwise,}
\end{array} \right.
\end{eqnarray}
where
$\xi_{01}(\underline{V})=\min\left\{\frac{1}{n_1},\frac{T_{2}}{T_{1}}\frac{1}{n_1+n_2}\right\}$
and
$\xi_{02}=\max\left\{\frac{1}{n_2},\frac{T_{1}}{T_{2}}\frac{1}{n_1+n_2}\right\}$
have uniformly smaller risk than the restricted maximum likelihood
estimators of $\sigma_1$ and $\sigma_2$, respectively.
\end{corollary}
\subsection{A numerical study}
Here we carry out a simulation study to compare the
risk functions of the proposed estimators numerically. We have generated $50,000$ random samples from two
independent exponential distributions using $R$ software. The
percentage relative risk improvement (PRRI) of an estimator
$\delta_{1}$ over another estimator $\delta_{2}$ is defined as
\begin{eqnarray}
PRRI(\delta_{1},\delta_{2})=\left(\frac{Risk(\delta_{2})-Risk(\delta_1)}{Risk(\delta_{2})}\right)
\times 100\%
\end{eqnarray}
In the following tables, we present PRRI of various estimators for
different sample sizes and different values of $p$. We have
considered sample sizes $(n_1,n_2)$ as $(3,5),$ $(10,15),$ $(8,5)$
and $(15,10).$ The values of $p$ are taken to be $0.5$ and $-0.5$.
The following points can be observed from Table $1$ and Table $2$.
\begin{enumerate}
\item With respect to the BSEE, the PRRI of the improved estimators
decreases when sample sizes $(n_{1},n_{2})$ increases with the
constraint $n_{1}<n_{2}$. Similar observation is noticed for
$n_{1}>n_{2}$ case.

\item For fixed $(n_{1},n_{2})$, the PRRI of the improved estimators (with respect to the
BSEE) increases if $(\sigma_{1},\sigma_{2})$ increases. Further, if
we increase the values of $\sigma_{2}$ by keeping $\sigma_{1}$
fixed, then the PRRI of the improved estimators (with respect to the BSEE)
decreases. This observation ensures that after some distance, the
risk function of the improved BSEE coincides with that of the BSEE.

\item It can seen that the risk performance of restricted MLE and improved restricted MLE  increases when the  sample sizes $(n_{1},n_{2})$ increases with the constraint $n_{1}<n_{2}$.   
\end{enumerate}
\begin{table}[ht]
\caption{The PRRI of various estimators for $k=2$ and $p=0.5.$ The
3rd, 4th and 5th columns respectively present the PRRI of improved
BAEE, restricted MLE and improved restricted MLE over the BAEE
whereas the 6th, 7th and 8th columns present the PRRI of improved
BAEE, BAEE and improved restricted MLE over the restricted MLE.}
\vspace{.3cm}
\centering 
\scalebox{0.82}{
\begin{tabular}{c c c c c c c c} 
\hline\hline 
$(n_{1},n_{2})$ & $(\sigma_{1},\sigma_{2})$ &
$(\delta_{\phi_{10}},\delta_{\phi_{20}})$ &
$(\delta^{2}_{1RML},\delta^{2}_{2RML})$ &
$(\delta^{2}_{1IRML},\delta^{2}_{2IRML})$ &
$(\delta_{\phi_{10}},\delta_{\phi_{20}})$ &
$(\delta_{1d_{01}},\delta_{2d_{02}})$
& $(\delta^{2}_{1IRML},\delta^{2}_{2IRML})$ \\ [0.5ex] 
\hline\hline\vspace{0.1cm} 
(3,5) & (0.2,0.5) & (4.35,6.54) & (12.26,-0.71) & (13.50,3.94) &
(-9.00,7.20)
& (-13.09,0.71) & (1.41,4.62) \\\vspace{0.1cm}
~ & (0.2,1.0) & (0.88,1.06) & (7.26,-1.01) & (7.44,-0.23) &
(-6.88,2.05)
& (-7.87,1.00) & (0.19,0.68) \\\vspace{0.1cm}
~ & (0.2,1.5) & (0.42,0.28) & (2.04,-1.03) & (2.44,-0.85) &
(-2.03,1.29)
& (-2.46,1.02) & (0.04,0.18) \\\vspace{0.1cm}
~ & (0.5,1.0) & (6.06,10.21) & (9.89,-0.37) & (12.22,7.01) &
(-4.25,10.55)
& (-10.98,0.37) & (2.58,7.35) \\\vspace{0.1cm}
~ & (0.5,1.5) & (3.01,4.30) & (12.04,-0.86) & (12.87,2.15) &
(-10.26,5.12)
& (-13.69,0.85) & (0.84,2.98) \\\vspace{0.1cm}
~ & (0.5,2.0) & (1.41,2.04) & (9.87,-0.97) & (10.12,0.39) &
(-9.27,2.98)
& (-10.87,0.96) & (0.38,1.35) \\\vspace{0.1cm}
~ & (1.0,1.5) & (8.26,15.98) & (1.69,0.65) & (6.53,12.11) &
(6.68,15.42)
& (-1.71,-0.66) & (4.92,11.53) \\\vspace{0.1cm}
~ & (1.0,2.0) & (6.06,10.21) & (9.89,-0.37) & (12.22,7.01) &
(-4.25,10.55)
& (-10.98,0.36) & (2.58,7.35) \\\vspace{0.1cm}
~ & (1.0,2.5) & (4.35,6.54) & (12.25,-0.71) & (13.50,3.94) &
(-9.00,7.20) & (-13.97,0.70) & (1.41,4.62) \\ \\\\
(10,15) & (0.2,0.5) & (0.83,0.90) & (15.79,-0.44) & (15.82,-0.24) &
(-17.76,1.35)
& (-18.76,0.44) & (0.03,0.69) \\\vspace{0.1cm}
~ & (0.2,1.0) & (0,0.002) & (0.18,-0.44) & (0.18,-0.45) &
(-0.18,0.44)
& (-0.18,0.44) & (0,0.001) \\\vspace{0.1cm}
~ & (0.2,1.5) & (0,0) & (-0.53,-0.44) & (-0.53,-0.45) & (0.52,0.45)
& (0.53,0.44) & (0,0) \\\vspace{0.1cm}
~ & (0.5,1.0) & (2.94,3.16) & (17.64,-0.44) & (17.85,2.11) &
(-17.84,3.59)
& (-21.42,0.44) & (0.25,2.54) \\\vspace{0.1cm}
~ & (0.5,1.5) & (0.15,0.26) & (9.59,-0.44) & (9.59,-0.25) &
(-10.43,0.70)
& (-10.61,0.44) & (0.01,0.19) \\\vspace{0.1cm}
~ & (0.5,2.0) & (0.01,0.02) & (2.42,-0.44) & (-2.42,-0.43) &
(-2.47,0.46)
& (-2.49,0.44) & (0,0.01) \\\vspace{0.1cm}
~ & (1.0,1.5) & (9.05,10.87) & (-6.79,-0.42) & (-5.15,8.79) &
(14.84,11.24)
& (6.36,0.42) & (1.54,9.16) \\\vspace{0.1cm}
~ & (1.0,2.0) & (2.94,3.16) & (17.64,-0.45) & (17.85,2.11) &
(-17.84,3.59)
& (-21.42,0.44) & (0.25,2.54) \\\vspace{0.1cm}
~ & (1.0,2.5) & (0.83,0.90) & (15.79,-0.45) & (15.82,0.25) &
(-17.76,1.34) & (-18.76,0.44) & (0.03,0.69) \\ \\\\
(8,5) & (0.2,0.5) & (2.73,8.67) & (-54.73,0.39) & (-7.14,5.99) &
(37.14,8.31)
& (35.37,-0.39) & (30.75,5.62) \\\vspace{0.1cm}
~ & (0.2,1.0) & (0.28,0.79) & (1.23,-1.07) & (8.78,-0.54) &
(-0.96,1.84)
& (-1.25,1.06) & (7.64,0.52) \\\vspace{0.1cm}
~ & (0.2,1.5) & (0.10,0.14) & (3.54,-1.13) & (5.49,-1.04) &
(-3.59,1.26)
& (-3.67,1.12) & (2.02,0.08) \\\vspace{0.1cm}
~ & (0.5,1.0) & (4.38,15.67) & (-95.99,2.40) & (-20.71,12.05) &
(51.21,13.59)
& (48.97,-2.46) & (38.41,9.88) \\\vspace{0.1cm}
~ & (0.5,1.5) & (1.68,5.04) & (-29.96,-0.42) & (1.26,2.88) &
(24.35,5.40)
& (23.06,0.42) & (24.03,3.29) \\\vspace{0.1cm}
~ & (0.5,2.0) & (0.62,1.87) & (-6.67,-0.94) & (8.00,0.31) &
(6.81,2.79)
& (6.22,0.93) & (13.37,1.24) \\\vspace{0.1cm}
~ & (1.0,1.5) & (6.96,28.68) & (-162.75,7.87) & (-39.43,23.68) &
(64.59,22.58)
& (61.94,-8.55) & (46.93,17.15) \\\vspace{0.1cm}
~ & (1.0,2.0) & (4.37,15.67) & (-95.99,2.40) & (-20.71,12.05) &
(51.21,13.59)
& (48.97,-2.46) & (38.41,9.88) \\\vspace{0.1cm}
~ & (1.0,2.5) & (2.73,8.67) & (-54.73,0.39) & (-7.14,5.99) &
(36.13,8.31) & (35.37,-0.39) & (30.75,5.62) \\ \\\\
(15,10) & (0.2,0.5) & (0.79,1.79) & (-35.54,-0.44) & (-13.89,0.95) &
(26.81,2.23)
& (26.22,0.44) & (15.97,1.39) \\\vspace{0.1cm}
~ & (0.2,1.0) & (0.005,0.01) & (5.55,-0.47) & (5.88,-0.47) &
(-5.87,0.48) & (-5.88,0.47) & (0.35,0.01) \\\vspace{0.1cm}
~ & (0.2,1.5) & (0,0) & (0.76,-0.47) & (0.79,-0.48) &
(-0.77,0.47) & (-0.77,0.47) & (0.30,0) \\\vspace{0.1cm}
~ & (0.5,1.0) & (2.25,5.80) & (-105.05,-0.24) & (-15.86,4.30) &
(52.32,6.03) & (51.23,0.25) & (25.94,4.53) \\\vspace{0.1cm}
~ & (0.5,1.5) & (0.27,0.59) & (-5.78,-0.46) & (3.08,-0.03) &
(5.72,1.05) & (5.47,0.56) & (8.38,0.44) \\\vspace{0.1cm}
~ & (0.5,2.0) & (0.03,0.10) & (7.78,-0.48) & (9.49,-0.44) &
(-8.40,0.54) & (-8.43,0.47) & (1.86,0.04) \\\vspace{0.1cm}
~ & (1.0,1.5) & (6.56,18.19) & (-242.13,1.25) & (-110.63,15.21) &
(72.69,17.15) & (17.77,-1.27) & (38.43,14.13) \\\vspace{0.1cm}
~ & (1.0,2.0) & (2.25,5.80) & (-105.05,-0.25) & (-51.86,4.30) &
(52.33,6.03) & (51.23,0.24) & (24.94,4.53) \\\vspace{0.1cm}
~ & (1.0,2.5) & (0.79,1.79) & (-35.54,-0.44) & (-13.89,0.95) &
(26.81,2.23) & (26.22,0.44) & (15.97,1.39) \\
[1ex] 
\hline\hline 
\end{tabular}}
\label{tb1} 
\end{table}
\clearpage

\begin{table}[ht]
\caption{The PRRI of various estimators for $k=2$ and $p=-0.5.$ The
3rd, 4th and 5th columns respectively present the PRRI of improved
BAEE, restricted MLE and improved restricted MLE over the BAEE
whereas the 6th, 7th and 8th columns present the PRRI of improved
BAEE, BAEE and improved restricted MLE over the restricted MLE.}
\vspace{0.4cm}
\centering 
\scalebox{0.82}{
\begin{tabular}{c c c c c c c c} 
\hline\hline \vspace{0.1cm}
$(n_{1},n_{2})$ & $(\sigma_{1},\sigma_{2})$ &
$(\delta_{\phi_{10}},\delta_{\phi_{20}})$ &
$(\delta^{2}_{1RML},\delta^{2}_{2RML})$ &
$(\delta^{2}_{1IRML},\delta^{2}_{2IRML})$ &
$(\delta_{\phi_{10}},\delta_{\phi_{20}})$ &
$(\delta_{1d_{01}},\delta_{2d_{02}})$
& $(\delta^{2}_{1IRML},\delta^{2}_{2IRML})$ \\ [0.5ex] 
\hline\hline\vspace{0.1cm} 
(3,5) & (0.2,0.5) & (5.18,6.98) & (3.01,-0.57) & (4.75,8.16) &
(2.25,7.50)
& (-3.09,0.56) & (1.81,8.67) \\ \vspace{0.1cm}
~ & (0.2,1.0) & (1.09,1.16) & (2.65,-0.98) & (2.89,0.71) &
(-1.61,2.11)
& (-2.72,0.96) & (0.25,1.67) \\\vspace{0.1cm}
~ & (0.2,1.5) & (0.44,0.31) & (0.45,-1.01) & (0.52,-0.51) &
(-0.02,1.31)
& (-0.46,0.99) & (0.06,0.48) \\\vspace{0.1cm}
~ & (0.5,1.0) & (7.32,10.77) & (-0.84,-0.11) & (2.44,12.65) &
(8.09,10.86)
& (0.84,0.11) & (3.25,12.74) \\\vspace{0.1cm}
~ & (0.5,1.5) & (3.69,4.64) & (4.07,-0.77) & (5.12,5.27) &
(-0.39,5.37)
& (-4.24,0.76) & (1.07,5.98) \\\vspace{0.1cm}
~ & (0.5,2.0) & (1.92,2.22) & (3.72,-0.92) & (4.19,2.15) &
(-1.87,3.12)
& (-3.87,0.91) & (0.48,3.05) \\\vspace{0.1cm}
~ & (1.0,1.5) & (10.14,16.63) & (-11.23,1.21) & (-4.45,19.45) &
(19.22,15.61)
& (10.10,-1.22) & (6.09,18.25) \\\vspace{0.1cm}
~ & (1.0,2.0) & (7.31,10.77) & (-0.84,-0.11) & (2.43,12.65) &
(8.09,10.87)
& (0.84,0.11) & (3.26,12.75) \\\vspace{0.1cm}
~ & (1.0,2.5) & (5.18,6.98) & (3.00,-0.57) & (4.75,8.16) &
(2.25,7.50) & (-3.09,0.56) & (1.81,8.68) \\ \\\\
(10,15) & (0.2,0.5) & (0.94,0.91) & (11.1,-0.35) & (11.1,0.83) &
(-17.76,1.34)
& (-18.76,0.44) & (0.03,0.69) \\\vspace{0.1cm}
~ & (0.2,1.0) & (0,0.002) & (-0.08,-0.08) & (-0.08,-0.34) &
(-0.18,0.45)
& (0.53,0.44) & (0,0) \\\vspace{0.1cm}
~ & (0.2,1.5) & (0,0) & (-0.56,-0.56) & (-0.56,-0.35) & (0.52,0.44)
& (0.53,0.44) & (0,0) \\\vspace{0.1cm}
~ & (0.5,1.0) & (3.23,3.17) & (10.45,-0.35) & (10.79,3.58) &
(-17.85,3.59)
& (-21.42,0.45) & (0.25,2.54) \\
~ & (0.5,1.5) & (0.23,0.26) & (6.71,-0.35) & (6.72,0.005) &
(-10.43,0.71)
& (-10.61,0.45) & (0.006,0.19) \\\vspace{0.1cm}
~ & (0.5,2.0) & (0.02,0.02) & (1.47,-0.35) & (1.47,-0.32) &
(-2.47,0.46)
& (-2.48,0.45) & (0,0.013) \\\vspace{0.1cm}
~ & (1.0,1.5) & (10.02,10.83) & (-19.06,-0.31) & (-16.81,12.4) &
(14.84,11.25)
& (6.36,0.41) & (1.54,9.17) \\\vspace{0.1cm}
~ & (1.0,2.0) & (3.22,3.17) & (10.49,-0.35) & (10.79,3.58) &
(-17.84,3.59)
& (-21.42,0.44) & (0.25,2.54) \\\vspace{0.1cm}
~ & (1.0,2.5) & (0.93,0.91) & (11.12,-0.35) & (11.14,0.83) &
(-17.76,1.35) & (-18.76,0.45) & (0.03,0.69) \\ \\\\
(8,5) & (0.2,0.5) & (3.09,8.9) & (-79.04,1.18) & (-15.82,9.99) &
(37.14,8.31)
& (35.37,-0.39) & (30.75,5.62) \\\vspace{0.1cm}
~ & (0.2,1.0) & (0.35,0.81) & (-6.05,-0.81) & (3.94,0.16) &
(-0.96,1.84)
& (-1.25,1.06) & (7.64,0.52) \\\vspace{0.1cm}
~ & (0.2,1.5) & (0.08,0.13) & (0.51,-0.89) & (3.09,-0.71) &
(-3.59,1.26)
& (-3.67,1.12) & (2.03,0.08) \\\vspace{0.1cm}
~ & (0.5,1.0) & (4.68,16.02) & (-130.67,3.78) & (-30.31,18.17) &
(51.21,13.59)
& (48.97,-2.46) & (38.41,9.88) \\\vspace{0.1cm}
~ & (0.5,1.5) & (1.92,5.18) & (-47.95,0.09) & (-6.5,5.53) &
(24.35,5.40)
& (23.06,0.41) & (24.03,3.29) \\
~ & (0.5,2.0) & (0.76,1.93) & (-17.64,-0.62) & (1.71,1.58) &
(6.81,2.79)
& (6.23,0.93) & (13.37,1.24) \\\vspace{0.1cm}
~ & (1.0,1.5) & (7.71,28.81) & (-215.37,10.58) & (-50.04,32.48) &
(64.59,22.58)
& (61.94,-8.55) & (46.93,17.15) \\
~ & (1.0,2.0) & (4.88,16.02) & (-130.68,3.78) & (-30.32,18.17) &
(51.21,13.59)
& (48.97,-2.46) & (38.41,9.88) \\\vspace{0.1cm}
~ & (1.0,2.5) & (3.09,8.95) & (-79.04,1.17) & (-15.82,9.98) &
(36.14,8.32) & (35.37,-0.39) & (30.75,5.62) \\ \\\\
(15,10) & (0.2,0.5) & (0.85,1.78) & (-51.33,-0.62) & (-23.48,1.56) &
(26.81,2.23)
& (26.22,0.44) & (15.97,1.39) \\\vspace{0.1cm}
~ & (0.2,1.0) & (0.008,0.01) & (3.73,-0.67) & (4.17,-0.66) &
(-5.87,0.49) & (-5.88,0.47) & (0.35,0.01) \\\vspace{0.1cm}
~ & (0.2,1.5) & (0,0) & (0.53,-0.67) & (0.56,-0.67) &
(-0.77,0.47) & (-0.77,0.47) & (0.30,0) \\\vspace{0.1cm}
~ & (0.5,1.0) & (2.45,5.72) & (-135.07,-0.36) & (-65.83,6.27) &
(52.32,6.03) & (51.23,0.25) & (25.94,4.54) \\\vspace{0.1cm}
~ & (0.5,1.5) & (0.29,0.57) & (-15.18,-0.65) & (-3.66,0.08) &
(5.72,1.05) & (5.47,0.56) & (8.39,0.44) \\\vspace{0.1cm}
~ & (0.5,2.0) & (0.03,0.06) & (3.75,-0.67) & (5.97,-0.59) &
(-8.41,0.54) & (-8.43,0.48) & (1.86,0.04) \\\vspace{0.1cm}
~ & (1.0,1.5) & (7.02,17.92) & (-303.68,1.55) & (-131.17,20.14) &
(72.69,17.15) & (17.77,-1.27) & (38.43,14.14) \\\vspace{0.1cm}
~ & (1.0,2.0) & (2.45,5.72) & (-135.07,-0.36) & (-65.82,6.27) &
(52.33,6.03) & (51.23,0.24) & (24.94,4.53) \\\vspace{0.1cm}
~ & (1.0,2.5) & (0.86,1.79) & (-51.54,-0.62) & (-23.48,1.56) &
(26.81,2.23) & (26.22,0.44) & (15.97,1.39) \\
[1ex] 
\hline\hline 
\end{tabular}}
\label{tb1} 
\end{table}
\clearpage

\section{Conclusion}
In the present communication, we have studied the component wise estimation of scale  parameters of several exponential distributions when the scale parameters are ordered with respect to the linex loss function. In Section \ref{sec2}, we have studied the estimation of ordered scale parameters when the location parameters are known. The inadmissibility of the usual estimator is proved by deriving a class of improved estimators. Using this result, we have obtained an estimator which improves upon the restricted MLE. When location parameters are unknown, the estimation of ordered scale parameters is investigated in section \ref{sec3}. It is proved that the BAEE is inadmissible by deriving a class of improved estimators. We have also established the inadmissibility of restricted MLE. Finally, for both cases, we have compared the risk performances of the proposed estimators through simulations.   
\bibliography{bibexplinex}
\end{document}